\documentclass[intlimits]{amsart}
\usepackage{amsmath}
\usepackage{amssymb}
\newcommand{\Sum}{\sum\limits}
\newtheorem{ksplemma}{Lemma}[section]
\newtheorem{ksprem}[ksplemma]{Remark}
\newtheorem{ksptheorem}[ksplemma]{Theorem}
\newtheorem{kspproposition}[ksplemma]{Proposition}
\newtheorem{kspdefinition}[ksplemma]{Definition}
\newtheorem{kspcorollary}[ksplemma]{Corollary}
\newcommand{\lcm}{\text{lcm}}
\newcommand{\Z}{\mathbb Z}
\newcommand{\Q}{\mathbb Q}
\newcommand{\C}{\mathbb C}
\numberwithin{equation}{section}
\def\mattwod(#1;#2){\left(\begin{array}{cc}
                               #1 & 0 \\
                               0  & #2
                                      \end{array}\right)}

\def\(#1;#2){\left(\begin{array} {c}                                 #1 \\                                #2 \end{array}\right)}

\def\mattwo(#1;#2;#3;#4){\left(\begin{array}{cc}
                               #1 & #2 \\
                               #3  & #4
                                      \end{array}\right)}

\begin{document}
\parindent=0pt
\title[Genus theta series, Hecke operators and the basis problem]{Genus theta series, Hecke operators and the basis problem for
  Eisenstein series}
\author{Hidenori Katsurada \and Rainer Schulze-Pillot}  
%\date{\today}
\thanks{Dedicated to the memory of Tsuneo Ararakawa}
\maketitle
\begin{abstract}{We derive explicit formulas for the action of the Hecke
  operator $T(p)$ on the genus theta series of a positive definite
  integral quadratic form and prove a theorem on the generation of
  spaces of Eisenstein series by genus theta series. We also discuss
  connections of our results with Kudla's matching principle for theta
integrals.}
\end{abstract}

\section{Introduction}
In the theory of theta series of positive definite quadratic forms the problem
of giving explicit formulas for the action of Hecke operators on theta
series has received some attention \cite{An,Yo}.

\vspace{0.3cm}
If $p$ is prime to the level $N$ of the quadratic form $q$ of rank $m$ in
question, the action of the usual generators $T(p), T_i(p^2)$ of the $p$-part
of the Hecke algebra for the group $\Gamma_0^{(n)}(N) \subseteq Sp_n
(\Bbb Z)$ is known \cite{An,Yo} except for the case that $n < \frac{m}{2}$ and
$\chi(p) = -1$, where $\chi$ is the nebentype character of the degree
$n$ theta series of $q$. In this last case it is unknown whether $T(p)$ 
leaves the space of cusp forms generated by the theta series of positive
definite quadratic forms of the same level and rational square class of
the discriminant invariant. Some deep results concerning this question
have been obtained by Waldspurger \cite{W}.

\vspace{0.3cm}
To our surprise, there seem to be no results available even for the question
how to describe the action of $T(p)$ on the genus theta series of $q$, i.e.,
Siegel's weighted average over the theta series of the quadratic forms
$q'$ in the genus of $q$.

\vspace{0.3cm}
The present note intends to fill this gap. It turns out that we have
different methods available to express the image of the genus 
theta series under the operator $T(p)$ in terms of theta series:
Using results of Freitag \cite{freitag-hamburg96}, Salvati Manni
\cite{manni-hamburg97} and Chiera \cite{chiera-hamburg03} one obtains an
expression as a linear combination of theta series of positive
definite quadratic forms of level $\lcm (N,4)$.

We show in Section 5 that this result can be improved to an (explicit)
expression as a linear combination of genus theta series of positive
definite quadratic 
forms of level $N$ if $N$ is an odd prime. 
In fact we prove in that case that any $n+1$ of the genera of
quadratic forms that are rationally equivalent to the given genus and
have level dividing $N$ yield a basis of the relevant space of
holomorphic Eisenstein series.
 
This can
be generalized to arbitrary square free level under a slightly
technical condition on the degree $n$ depending on the
$\Q_p$-equivalence class for $p$ dividing $N$ of the given genus of quadratic
forms; generalizations to arbitrary level will be the subject of
future work.

\medskip
On the other hand, using the explicit expression for the action of
Hecke operators on Fourier coefficients of modular forms given in
\cite{An}, Siegel's mass formula and relations between the
local densities of quadratic forms we find a much simpler expression:
The genus theta series is transformed into a multiple of the genus theta series
of a different genus of quadratic forms. If $\chi(p)=-1$, the genus
involved turns out to be indefinite, and the theta 
series is the one defined by Siegel ($n=1$) and Maa\ss\
\cite{Si,Ma}. This phenomenon is an instance (with quite explicit data)
of the matching 
principle for Siegel-Weil integrals attached to different quadratic
spaces that has been observed by Kudla in \cite{kudla-compositio03},
we discuss this in Section 6.

\medskip As a consequence of our work we are able to give a positive
solution to the basis problem for modular forms in a number of new
cases; this will be done in joint work with S. B\"ocherer.

 \section{Preliminaries}
Let $L$ be a lattice of full rank on the $m$-dimensional vector space
$V$ over $\Bbb Q$, $q:\: V \longrightarrow \Bbb Q$ a positive definite
quadratic form with $q(L) \subseteq \Bbb Z$, $B(x,y) = q(x+y)-q(x)-q(y)$
the associated symmetric bilinear form, $N=N(L)$ the level of $q$ (i.e., 
$N^{-1}\Z = q(L^{\#})\Z$, where $L^{\#}$ is the dual lattice of
$L$ with respect to $B$); we assume $m = 2k$ to be even. 

Let $R$ be $\Z$ or $\Z_p$ for some prime $p$ and let ${\mathcal
  H}_n(R)$ denote the set of half-integral matrices of 
degree $n$ over $R,$ that is, ${\mathcal H}_n(R)$ is the set of
symmetric matrices $(a_{ij})$ of degree $n$  with entries in $\frac{1}{2}R$
such that $a_{ii} \ (i=1,...,n)$ and $2a_{ij} \ (1 \le i \not= j \le
n)$ belong to $R.$  

We note that for ${\bf x}=(x_1,\dots,x_n)  \in L^n$ the matrix 
$q({\bf x}):=
(\frac{1}{2} B(x_i,x_j))$ is in the set ${\mathcal H}_n(\Z)$; we also
note that ${\mathcal 
  H}_n(\Z_p)$ is equal to the 
set $M_n^{\rm sym}(\Z_p)$ of symmetric $n\times n$ matrices over $\Z_p$ for $p
\ne 2.$ For two square matrices $T_1$ and $T_2$ we write $T_1 \perp T_2
=\left(\begin{smallmatrix}
T_1&0\\0&T_2\end{smallmatrix}\right).$ 

We often write $a \perp T$ instead of $(a) \perp T$ if
$(a)$ is  a matrix of degree 1. If $K=(K,q')$ is a quadratic
$\Z_p$-lattice with Gram matrix 
$T$ with respect to some basis we will freely switch notation between
$T$ and $K$, so for example if $K$ is a one-dimensional lattice with
basis vector of squared length $a$ and $M$ a quadratic lattice with Gram
matrix $T$ we write as above  $a \perp T=(a) \perp T=K \perp T =
K\perp M.$

The theta series
 $$\vartheta^{(n)}(L,Z) = \sum_{{\bf x}=(x_1,\ldots,x_n)\in L^n}
  \exp(2\pi i~{\rm tr}(q({\bf x}) Z)$$
of degree $n$ of $(L,q)$ is well-known to be in the space 
$M_k^{(n)}(\Gamma_0^{(n)}(N),\chi)$ of Siegel modular forms of weight
$k = \frac{m}{2}$ and character $\chi$, where $\chi$ is the character of
$\Gamma_0^{(n)}(N)$ given by $\chi\left({A\,B\choose C\,D}\right) =
\tilde{\chi}(\det D)$, $\tilde{\chi}$ is the Dirichlet character modulo
$N$ given by $\tilde{\chi}(d) = \left(\frac{(-1)^k\det L}{d}\right)$ for
$d > 0$ and $\det L$ is the determinant of the Gram matrix of $L$ with
respect to some basis \cite{An}. 

For definitions and notations
concerning modular forms we refer again to \cite{An}, we recall that  the Hecke
operator associated to the double coset 
\begin{equation*}
\Gamma_0(N) \begin{pmatrix} 1 & &&&&\\
&\ddots&&&&\\
&&1&&&\\
&&&p&&\\
&&&&&\ddots&\\
&&&&&&p
\end{pmatrix} \Gamma_0(N)
\end{equation*} 
is as usual denoted by $T(p).$
 \vspace{0.3cm}\\
We let $\{L_1,\ldots,L_h\}$ be a set of representatives of the classes
of lattices in the genus of $L$, put $w = \sum_{i=1}^{h} 
\frac{1}{|O(L_i)|}$ (where $O(L_i)$ is the group of isometries of $L$
onto itself with respect to $q$) and write
 $$\vartheta^{(n)}({\rm gen}~L,Z) = \frac{1}{w} \sum_{i=1}^{h}
  \frac{\vartheta^{(n)}(L_i,Z)}{|O(L_i)|}$$
for Siegel's weighted average over the genus.
 \vspace{0.3cm}\\
By Siegel's theorem (see \cite{Ki})
the Fourier coefficient $r({\rm gen}~L,A)$ at a
positive semidefinite half integral symmetric matrix $A$ can be
expressed as a product of local densities,
 \begin{equation}\label{siegel_theorem}
 r({\rm gen}~L,A) = c\cdot(\det A)^{\frac{m-n-1}{2}} (\det L)^{\frac
 {n}{2}} \prod_{\ell~{\rm prime}} \alpha_{\ell} (L,A)
 \end{equation}
with some constant $c.$

Here the local density $\alpha_{\ell} (L,A)$ is given as 
 \begin{eqnarray*}
 \alpha_{\ell}(L,A) &= &\alpha_{\ell}(S,A)\\
&=&\ell^{j\cdot(\frac{n\cdot(n+1)}{2}-mn)} \cdot
  \#\{{\bf x}\in L^n/\ell^jL^n~|~q({\bf x}) \equiv A\bmod \ell^j
  {\mathcal H}_n(\Z_{\ell} \}\\
&=& \ell^{j\cdot(\frac{n\cdot(n+1)}{2}-mn)} \# {\mathcal A}_j(S,A),
 \end{eqnarray*}
for sufficiently large $j$ with  an additional factor $\frac{1}{2}$ if
$m = n$ where $S$ denotes a Gram matrix of $L$ and
where we write
\begin{eqnarray*}
{\mathcal A}_j(L,A)&=&{\mathcal A}_j(S,A)\\
&=&\{{\bf x}\in L^n/\ell^jL^n~|~q({\bf x}) \equiv A\bmod \ell^j
  {\mathcal H}_n(\Z_{\ell} \}\\
&=&\{ X =(x_{ij}) \in
M_{m,n}({\Z}_{\ell})/\ell^jM_{m,n}({\Z}_{\ell})\mid A[X]- B \in
  \ell^j{\mathcal H}_n({\Z}_{\ell}) \} 
\end{eqnarray*}

\section{Eisenstein series and theta series}
\begin{kspproposition}\label{freitag-manni-proposition}
Let $L$ be a lattice of rank $m=2k$ with positive definite quadratic
form $q$ of square free level $N$, let $n<k-1$ and let
$F=\vartheta^{(n)}(\text{gen}(L))$ denote the genus theta series of
$L$ of degree
$n$. Then for any prime $p\nmid N$ the modular form $F\vert_kT(p)$ is a
linear combination of genus theta series of genera of lattices with
positive definite quadratic form of level $N'=\lcm(N,4).$
\end{kspproposition}
\begin{proof} By \cite{An2} $G:=F\vert_kT(p)$ is an eigenfunction of
  infinitely many Hecke operators $T(\ell)$ for the primes $\ell \nmid
  pN$ with $\chi(\ell)=1$ (where $\chi$ is the nebentyp character
  for $\vartheta^{(n)}(L)$). Proposition 4.3 of
  \cite{freitag-hamburg96} implies then that $G$ is in the
  space that is generated by Eisenstein series for the
  principal congruence subgroup of level $N$; this can also
  be obtained from Siegel's main theorem if one uses that this space
  is Hecke invariant. 

We want now to use 
  Theorem 6.9 of
  \cite{freitag-hamburg96} (see also \cite{manni-hamburg97}) to
  prove that $G$ is a linear combination of theta series with
  characteristic for the principal congruence subgroup of level
  $N'=\lcm(N,4)$. For this recall that with 
\begin{eqnarray*}
\Gamma_1^{(n)}(N')&:=&\{\begin{pmatrix}
A&B\\C&D
\end{pmatrix} \in Sp_n(\Z) \mid C \equiv 0 \bmod N', \det(A) \equiv 1
\bmod N'\}\\ 
\Gamma_{sq}^{(n)}(N')&:=&\{\begin{pmatrix}
A&B\\C&D
\end{pmatrix} \in Sp_n(\Z) \mid C \equiv 0 \bmod N',\\ &&\quad\quad \det(A) \text{ is
  congruent to a square }
\bmod N'\}\\
\Delta_{sq}^{(n)}(N')&:=&\{\begin{pmatrix}
A&B\\C&D
\end{pmatrix} \in Sp_n(\Z) \mid C \equiv B \equiv 0 \bmod N',\\
&&\quad\quad \det(A) \text{ is
  congruent to a square }
\bmod N'\}
\end{eqnarray*}
we have the cusp function
$\Phi_G:Sp_n(\Z)/
\Gamma_1^{(n)}(N') \longrightarrow \C$
associating to the class $\gamma \Gamma_1^{(n)}(N')$ the value of $G$
in the cusp $\gamma \infty$. The quoted theorem then states that $G$
is a
linear combination of theta series as above if and only if $\Phi_G$ is
constant on all cosets $\gamma\Gamma_{sq}^{(n)}(N')$.
Since $N$ is square free it is known that we can choose a set of
representatives $w_i\infty$ of the classes of
$\Gamma_0^{(n)}(N)$-equivalent cusps where the $w_i$ are certain
involutions normalizing the group $\Delta_{sq}^{(n)}(N)$, see e.g. 
\cite[Lemma 8.1]{boe-sp_nagoya}.
Since we have
$\Gamma_{sq}^{(n)}(N)=\Delta_{sq}^{(n)}(N)\Gamma_1^{(n)}(N)$
and since $G$ transforms under $\Gamma_0^{(n)}(N)$ according to the
quadratic character $\chi$ one sees that $\Phi_G$ is constant on the
cosets $\gamma\Gamma_{sq}^{(n)}(N)$ and hence also on the cosets $\gamma\Gamma_{sq}^{(n)}(N')$.

Our modular form $G$ is therefore indeed a linear combination of theta series with
  characteristic for the principal congruence subgroup of level
  $N'=\lcm(N,4)$.
Since $G$ is in fact a modular form for
  $\Gamma_0(N'),$  Chiera's Theorem 1 \cite{chiera-hamburg03} implies
  that $G$ is a linear combination of theta series
  $\vartheta^{(n)}(K_j)$ attached to full lattices $K_j$ with
  quadratic form of level dividing $N'$.

 It is well known
  that the values of the theta series of lattices in the same genus at
  zero dimensional cusps are the same. From Proposition 3.3 of
  \cite{freitag-hamburg96} we can then conclude that $G$ is in fact a
  linear combination of the $\vartheta^{(n)}(\text{gen}(K_j))$ as
  asserted.
\end{proof}

\section{Action of $T(p)$ and local densities}
The action of the Hecke operator $T(p)$ on the Fourier
coefficients of a Siegel modular form at nondegenerate matrices $A$
has been described explicitly 
by Maa\ss\ \cite{maass-ma-1951} and by Andrianov (Ex. 4. 2. 10 of \cite{An}):
\vspace{0.3cm}\\
Let $K$ be a $\Bbb Z$-lattice with quadratic form of rank $n$ that has
Gram matrix $p \cdot A$ with respect to some basis and write ${\mathcal M}_i$
for the set of lattices $M \supset K$ for which $K$ has elementary divisors
$(1,\ldots,1,p,\ldots,p)$ with $(n-i)$ entries $p$.\\
Then if $F(Z) \in M_k^{(n)}(\Gamma_0^{(n)}(N),\chi)$,
 $$\begin{array}{lrl}
  F(Z) & = & \sum_{A\ge 0} f(A) \exp (2\pi i {\rm tr}(AZ)),\vspace{0.2cm}\\
  G(Z) & = & (F|_k T(p))(Z) = \sum g_p(A) \exp(2 \pi i {\rm tr}(AZ)),
 \end{array}$$
one has for non-degenerate $A$:
 \begin{equation}\label{andrianov_formula}
  g_p(A) = \chi(p)^n p^{nk-\frac{n(n+1)}{2}} \sum_{i=0}^{n}(\bar{\chi}(p)
  p^{-k})^i p^{i \frac{i+1}{2}} \sum_{M\in {\mathcal M}_i} f(M),
 \end{equation}
where by $f(M)$ we denote the Fourier coefficient at an arbitrary Gram matrix
of the lattice $M$ (the coefficient $f(A)$ depends only on the integral 
equivalence class of $A$). Here by convention $f(M)$ is zero if the Gram
matrix of $M$ is not half integral.

\begin{kspproposition} Let 
\begin{equation*}F(Z) := \vartheta^{(n)}({\rm gen}~L,Z) =
   \sum_{A \geq 0} f(A) \exp 
(2\pi i~{\rm tr}(AZ)),
\end{equation*}
 $$G(Z) := (F|_k T(p)) (Z) = \sum_{A\geq 0}g_p(A) \exp(2\pi i~{\rm tr}(AZ)).$$
Then $g_p(A) = \lambda_p(L) (c \cdot (\det L)^{\frac{n}{2}} (\det
A)^{\frac{m-n-1}{2}}\prod_{\ell~{\rm prime}} \alpha_{\ell}  
(\tilde{L}_{\ell},A))$, where
 $$\lambda_p(L) = p^{nk-\frac{n(n+1)}{2}}\prod_{j=1}^{n} (1+\chi(p)p^{j-k})$$
and the $\Bbb Z_p$-lattice $\tilde{L}_{\ell}$ is given by
 $$\tilde{L}_{\ell} = \left\{ \begin{array}{rl}
  L_{\ell} & \mbox{ if } p = \ell\\
  ^pL_{\ell} & \mbox{ otherwise}. \end{array}\right.$$
Here $^pL_{\ell}$ denotes the lattice $L_{\ell}$ with quadratic form scaled
by $p$.
\end{kspproposition}
\begin{proof}
It is (by induction) enough to consider nondegenerate $A.$
We write the total factor in front of $f(M)$ for $M \in {\mathcal M}_i$
in (\ref{andrianov_formula})  as $\gamma_i$ and rewrite (\ref{andrianov_formula})  in the present situation as
 \begin{equation}
 g_p(A) = c\cdot (\det L)^{\frac{n}{2}} \sum_{i=0}^{n}
 \gamma_i
 \sum_{M\in {\mathcal M}_i} (\det M)^{\frac{m-n-1}{2}}
 \prod_{\ell} \alpha_{\ell} (L_{\ell},M)
 \end{equation}
by inserting the expression for $f(M)$ from (\ref{siegel_theorem})
(Siegel's theorem).
 \vspace{0.3cm}\\
Since $\det M = p^{2i-n} \det A$ for $M \in {\mathcal M}_i$ this becomes
 \begin{equation}\begin{split}
 g_p(A) = &c \cdot (\det L)^{\frac{n}{2}} (\det A)^{\frac{m-n-1}{2}} 
 p^{-n({\frac{m-n-1}{2}})}\\&\cdot \sum_{i=0}^{n} \gamma_i 
\sum_{M \in {\mathcal M}_i} p^{(\frac{m-n-1}{2}) \cdot 2i}
 \prod_{\ell} \alpha_{\ell} (L_{\ell},M_{\ell}).
 \end{split}\end{equation}
Now for $\ell \not= p$ we have $M_\ell = K_\ell$ for all $M$ occurring, hence
$\alpha_{\ell}(L_{\ell},M_{\ell}) = \alpha_{\ell}(L_{\ell},pA) 
= \alpha_{\ell}(^pL_{\ell},A) = \alpha_{\ell}(\tilde{L}_{\ell},A)$,
for all $\ell \not= p$.
 \vspace{0.3cm}\\
So it remains to prove
 \begin{equation}\label{density_formula_1}
p^{-n\frac{m-n-1}{2}} \sum_{i=0}^{n} \gamma_i \sum_{M \in {\mathcal M}_i}
 p^{(\frac{m-n-1}{2}) \cdot 2i} \alpha_p(L_p,M_p) = \lambda_p(L)
 \alpha_p(L_p,A).
 \end{equation}

We insert
$\gamma_i=\chi(p)^np^{nk-\frac{n(n+1)}{2}}(\chi(p)p^{-k})^{i}p^{i\frac{i+1}{2}}$ 
and divide both sides of (\ref{density_formula_1}) by
$p^{nk-\frac{n(n+1)}{2}}$ to get
\begin{equation}\label{density_formula_2}
\begin{split}p^{\frac{n(n+1)}{2}}\sum_{i=1}^n(\chi(p)p^{-k})^{n-i}&
p^{-i(n+1)}p^{\frac{i(i+1)}{2}}\sum_{M \in {\mathcal M}_i}\alpha_p(L,M)\\
& =\prod_{j=1}^n(1+\chi(p)p^{j-k})\alpha_p(L,A) 
\end{split}
\end{equation}
as the assertion that we have to prove.

For $\chi(p)=1$ this is proved in \cite{Yo} (see also \cite{An2}),
where it is also proved for $\chi(p) = -1$ and $n \geq k$ (in which case
the factor $\lambda_p(L)$ is zero).
To prove it for $\chi(p) = -1$ notice that $L_p$ is unimodular even by
assumption.
By Lemma 3.5 of \cite{Shi} there exists a polynomial $G_p(M;X)$ such that 
$\alpha_p(\hat{L}_p,M) = G_p(M;\chi_{\hat{L}_p}(p) p^{-\hat{k}})$
is true for all (even) unimodular $\Bbb Z_p$-lattices $\hat{L}_p$
of even rank $2\hat{k}$ with $\hat{k} \in \Bbb N$ and with 
 $$\chi_{\hat{L}_p}(p) := \left\{ \begin{array}{rl}
  1 & \mbox{ if $(-1)^{\hat{k}} \det \hat{L}_p$ is a square in 
      $\Bbb Q_p$}\\
 -1 & \mbox{ otherwise.} \end{array}\right. $$
Hence both sides of our assertion (\ref{density_formula_2}) are polynomials in
$X = \chi(p)p^{-\hat{k}}$ as $\hat{L}_p$ varies over (even) unimodular
$\Bbb Z_p$-lattices of (varying) rank $2\hat{k}$. The truth of the 
assertion for $\hat{L}_p$ with $\chi_{\hat{L}_p}(p) = 1$ and 
$\hat{k}$ arbitrary shows that these polynomials take the same value
at infinitely many places, hence must be identical. The assertion is therefore
true for all even unimodular $L_p$ of even rank.
\end{proof}
 \begin{ksplemma}\label{spacelemma}
 There is a unique isometry class of rational quadratic spaces 
$\tilde{V} = (\tilde{V},\tilde{q})$ of dimension $m$, such that
 \begin{equation}\label{eq:tildeV}
 \tilde{V}_{\ell} \cong V'_\ell:=\left\{ \begin{array}{rl}
  ^pV_{\ell} & \mbox{ if } p \not= \ell \\
  V_p & \mbox{ if } p = \ell
 \end{array}\right. 
 \end{equation}
for finite primes $\ell$ and $\tilde{V}_{\infty} = \tilde{V} \otimes_{\Bbb Q}
\Bbb R$ is either positive definite or of signature $(m-2,2)$.\\
$\tilde{V}$ carries a lattice $\tilde{L}$ such that
 \begin{equation}\label{eq:tildeL}
 \tilde{L}_{\ell} \cong \left\{ \begin{array}{rl}
  ^pL_{\ell} & \mbox{ if } p \not= \ell \\
  L_p & \mbox{ if } p = \ell . \end{array}\right.
 \end{equation}
$\tilde{V}_{\infty}$ is indefinite if and only if $\chi(p) = -1$.
The same assertion is true if one requires $\tilde{V}_\infty$ to be of
signature $(m-2-4j,2+4j)$ instead of $(m-2,2)$ for some $1\le j \le
\frac{m-2}{4}.$
 \end{ksplemma}
{\it Proof.} If $s_{\ell} V_{\ell}$ denotes the Hasse symbol of the quadratic
space $V_{\ell}$ and $V'_{\ell}$ is the quadratic $\Bbb Q_\ell$-space
as in (\ref{eq:tildeV}),
the discriminant of $V'_{\ell}$ is that of $V_{\ell}$ and the product of
the Hasse symbols $s_{\ell} V'_{\ell}$ over the finite primes $\ell$ is
the Hilbert symbol
 $$(p,(-1)^{\frac{m}{2}} \det L)_p \cdot \prod_{\ell~{\rm prime}} s_{\ell}
  V_{\ell}$$
by Hilbert's reciprocity law, with $(p,(-1)^{\frac{m}{2}} \det L)_p = \chi(p)$.
 \vspace{0.3cm}\\
If $V'_{\infty}$ is positive definite for $\chi(p) = 1$ and of signature
$(m-2,2)$ if $\chi(p) = -1$ one sees therefore that 
${\rm disc}~V'_{\ell} = {\rm disc}~V_{\ell}$ for all $\ell$ (including
$\infty$) and $\prod_{\ell,\infty} s_{\ell} V'_{\ell} = 1$, hence there is
a rational quadratic space $\tilde{V}$ such that $\tilde{V}_{\ell} \cong
V'_{\ell}$ for all $\ell$ including $\infty$. The uniqueness of $\tilde{V}$ 
is clear from the Hasse-Minkowski theorem, and that $\tilde{L}$ as in
(\ref{eq:tildeL}) exists on $\tilde{V}$ is obvious.

\vspace{0.3cm}
We recall that for an integral lattice of positive determinant and even rank
Siegel \cite{Si} for degree one and Maa\ss\ \cite{Ma} for arbitrary degree
defined a holomorphic theta series in the indefinite case whose Fourier 
coefficients at positive definite $A$ 
are proportional to the product
of the local densities of 
that lattice, subject to the restriction that the signature
  $(m_+,m_-)$ satisfies the condition
  $\min(\frac{m_++m_--3}{2},m_+,m_-)\ge n$. 
Denote this theta series (if it is defined) for $\tilde{L}$, normalized such that
its Fourier coefficient at $A$ is equal to 
 $$r({\rm gen}~\tilde{L},A):=c \cdot (\det A)^{\frac{m-n-1}{2}} (\det \tilde{L})^{\frac{n}{2}}
  \prod_{\ell~{\rm prime}} \alpha_{\ell} (\tilde{L},A),$$
by $\vartheta(\tilde{L},Z)$ or also by $\vartheta({\rm gen}~\tilde{L}),Z)$
  (notice that this theta series does indeed depend only on the genus
  of the lattice). The signature condition is in our
  situation always satisfied if $n=1,$ for bigger $n$ it can be
  satisfied by choosing $j$ in \ref{spacelemma} appropriately if $n\le
  k-2 $ (with $k=m/2).$ If the signature condition is not
  satisfied, we use the same notation $r({\rm gen}~\tilde{L},A)$
  (without knowing a priori whether these numbers are 
  the Fourier coefficients of a modular form).
 \vspace{0.3cm}\\
Then we arrive at the following final result:
 \begin{ksptheorem} \label{heckeaction_indefinite} 
Let $L$ be as above, $p$ a prime with $p {\not|} \det L$, $\tilde{L}$
a quadratic lattice with
 $$\tilde{L}_{\ell} = \left\{ \begin{array}{rl}
  ^pL_{\ell} & \mbox{ if } p \not= \ell\\
   L_p & \mbox{ if } p = \ell. \end{array}\right. $$
and of signature $(m,0)$ if $\chi(p)=1$, of signature $(m-2,2)$ if
$\chi(p) = -1$. Then 
 $$\vartheta^{(n)} ({\rm gen}~L)~|~T(p) =\lambda_p(L) \vartheta^{(n)} ({\rm gen}~\tilde{L})$$
with 
 $$\lambda_p(L) =p^{nk-\frac{n(n+1)}{2}}
 \prod_{j=1}^n(1+\chi(p)p^{j-k}),$$
where $\vartheta^{(n)} ({\rm gen}~
  \tilde{L},Z)$ is a
 holomorphic modular form of the same level as $L$ 
whose Fourier coefficient at a positive definite
  matrix $A$ is equal to $r({\rm gen}~\tilde{L},A)$.
The modular form $\vartheta^{(n)} ({\rm gen}~
  \tilde{L},Z)$ is the usual genus theta series if $\tilde{L}$ is
 positive definite and is equal to the theta series of Siegel and
 Maa\ss\ from above if  $\tilde{L}$ is indefinite and this series is
 defined.  

In particular, for
  all  $n<k$ there exists a holomorphic modular form of the same level
 as $L$ with Fourier coefficients $r({\rm gen}~\tilde{L},A)$ at at positive definite
  matrices $A$.
\end{ksptheorem}
{\bf Remark.}
a) $\lambda_p(L) =0$ if $n\ge k$ holds with $\chi(p)=-1$, which agrees
with Andrianov's result \cite{An2} for this case.
\vskip0.3cm
b) In the introduction we mentioned the question whether the space of
cusp forms generated by the theta series of positive definite lattices of
fixed level and rational square class of the discriminant is invariant under 
the action of the Hecke operators. In view of our theorem we might 
reformulate this question by substituting ``modular forms'' for
``cusp forms'' and omitting the restriction to positive definite lattices.
Since the indefinite theta series of Siegel and Maa\ss\ don't 
contribute to the space of cusp forms, this doesn't change the problem with 
regard to the subspace of cusp forms.
 \vspace{0.3cm}\\
c) Of course the same result holds true when we take an indefinite lattice
$\tilde{L}$ of signature $(m-2,2)$ as above as our starting point. The 
lattices appearing in $\vartheta({\rm gen}~\tilde{L},z)~|~T(p)$ are then
positive definite if $\chi(p) = -1$, indefinite if $\chi(p) = +1$.

\section{Spaces of genus theta series for odd prime level}

We will need some additional notations in this section.\\
Let $p$ be an odd  prime.
 For a non-zero element $a \in {\Q}_p$ we put $\chi_p(a)=1,-1,$ or
 0 according as ${\Q}_p(a^{1/2})={\Q}_p, {\Q}_p(a^{1/2})$ is
 an unramified quadratic extension of ${\Q}_p,$ or ${\Q}_p(a^{1/2})$
 is a  ramified quadratic extension of ${\Q}_p.$
 For a non-degenerate half-integral matrix $B$ of even degree $n,$ put
 $\xi_p(B)=\chi_p((-1)^{n/2} \det B).$ 

Further for non-negative integers $l,e$ and matrices $A\in {\mathcal
  H}_m(\Z_p), B \in {\mathcal H}_n(\Z_p)$ define 
$${\mathcal B}_e(A,B)^{(l)} = \{ X =(x_{ij}) \in {\mathcal A}_e(A,B); \ {\rm
  rank}_{{\Z}_p/p{\Z}_p}( x_{i,j})_{1 \le i \le m,1 \le j \le
  l}= l \}$$ 
(with $ {\mathcal A}_e(A,B)$ as in Section 2) and 
$$ \beta_p(A,B)^{(l)} = \lim_{e \rightarrow \infty}
p^{(-mn+n(n+1)/2)e}\#{\mathcal B}_e(A,B)^{(l)}.$$ 
We note that 
$$\beta_p(A,B)^{(0)}=\alpha_p(A,B).$$
In particular put  
$$\beta_p(A,B)=\beta_p(A,B)^{(n)},$$
and call it (as usual) the primitive density.  Further for $0 \le i \le m$ put 
$$\pi_{m,i}= GL_m({\Z}_p) (pE_i \perp E_{m-i}) GL_m({\Z}_p)$$ 
Furthermore let $H_k = \overbrace{H \perp ...\perp H}^k$ with $H =\mattwo(0;1/2;1/2;0).$

\bigskip

 Our goal in this section is to prove the following theorem:\\
\begin{ksptheorem}\label{maintheorem}
% Let $N$ be an odd square free positive integer and let $(V,q)$ be a
% nondegenerate quadratic space over ${\bf Q}$ of even dimension $m=2k$ and
% of positive discriminant such that the anisotropic kernel of the $p$-adic 
% completion $V_p$ has at most  dimension $2$ for all $p\mid N$ and such
% that $V$ carries an even integral lattice of level dividing $N$. 

% Then for  $n<k$ the space of modular forms for
% $\Gamma_0^{(n)}(N)$ spanned by the genus theta series of degree $n$
% attached to the genera of lattices of level dividing $N$ on $V$ span a
% space of dimension $n+1$ which is equal to the space of Eisenstein
% series for the group 
% $\Gamma_0^{(n)}(N)$ of weight $k$ and quadratic character attached to
% $\text{disc }(V)$.

Let $p$ be an odd prime, $k,n \in {\bf N}$ with $n\le k-1$  and $p \equiv
(-1)^{k} \bmod 4.$ 

Then the space of modular forms for
$\Gamma_0^{(n)}(p)$ spanned by the genus theta series  of degree $n$
attached to the genus of positive
definite integral quadratic lattices of rank $2k$, level $p$ and discriminant
$p^{2r+1}$ for some  $0\le r<k$ and the space spanned by the genus
theta series of degree $n$ (in the sense of
Theorem \ref{heckeaction_indefinite}) attached to the genus of
integral quadratic lattices of signature $(2k-2-4j,2+4j)$ (with $1\le j \le
\frac{2k-2}{4}$ fixed), level $p$ and discriminant  
$p^{2r+1}$  for some $0\le r<k$ coincide. This space has dimension $n+1$ and
is equal to the space of holomorphic Eisenstein series for the group
$\Gamma_0^{(n)}(p)$ of weight $k$ and nontrivial quadratic character.

For each of these  signatures  the theta series of any $n+1$ of the
$k$ genera of level dividing $p$ and having this  signature form a
basis of this space of modular forms.

% Consequently, for any positive definite lattice $L$ on $V$ of level
% dividing $N$ as above and
% $F=\vartheta^{(n)}(\text{gen}(L))$  the modular form $F\vert_kT(\ell)$
% can be expressed as a linear combination of theta series of positive
% definite lattices of level dividing $N$ on $V$ for all primes $\ell \nmid N.$

\end{ksptheorem}
The proof of this theorem will require a few intermediate results 
which may be of independent interest. A half-integral matrix $S_0$ over
${\Z}_p$ is called  ${\Z}_p$-maximal if it is the empty matrix or a
matrix corresponding to a ${\Z}_p$-maximal lattice. The main result we
need is the following theorem, whose proof again is broken up into
several steps: 

\begin{ksptheorem}\label{densitypolynomial}
 Let $p$ be an odd prime, let $T \in {\mathcal H}_n({\Z}_p).$ Let $k$
 be a positive integer, and 
 $S_0$ be a ${\Z}_p$-maximal half-integral matrix of degree  not greater
 than $2.$ Then there exist rational numbers 
 $a_i=a_i(k,S_0,T) \ (i=0,1,2,...,n)$ such that  
$$\alpha_p(H_{k-l-1} \perp pH_l \perp S_0,T)=a_0+a_1p^l+...+a_np^{nl}$$
for any $l=0,1,...,k-1.$
\end{ksptheorem}
% {\bf -----------------------------\\ ?? Maybe this can be stated in
% greater generality: The nicest 
%   statement would be:\\
% Let $K$ be any $\Z_p$-maximal lattice of rank $1$ or $2.$\\
% Then there exist rational numbers
%  $a_i=a_i(k,K,T) \ (i=0,1,2,...,n)$ such that  
% $$\alpha_p(H_{k-l-1} \perp pH_l \perp K,T)=a_0+a_1p^l+...+a_np^{nl}$$
% for any $l=0,1,...,k-1.$??\\---------------------------------}

\bigskip

To prove the theorem, first we remark that for $p \ne 2$ a
${\Z}_p$-maximal matrix $S_0$ of degree not greater than $2$ is
equivalent over ${\Z}_p$  to one of the following matrices: 

(M-1) $\phi$ (empty matrix),

(M-2) $u_1$  with $u_1\in {\Z}_p^*,$  

(M-3) $pu_1$  with $u_1\in {\Z}_p^*,$  

(M-4)  $u_1 \perp  u_2$ with $u_1, u_2 \in {\Z}_p^*,$ 

(M-5)  $u_1 \perp pu_2$   with $u_1, u_2 \in {\Z}_p^*,$

(M-6)   $pu_1 \perp pu_2$   with $u_1, u_2 \in {\Z}_p^*$ such that $-u_1u_2 \not\in ({\Z}_p^*)^2$

\bigskip

\begin{ksplemma}\label{lemma1}
Let $S_0$ be the matrix in Theorem 5.2. For a non-negative integer $l$
put $B_l=B_{S_0,l}=pH_l \perp S_0$ and $\tilde B_{l,i}= \tilde B_{S_0,l,i}=H_i \perp
pH_{l-i} \perp S_0.$  Let $T \in {\mathcal H}_n({\Z}_p) \cap
GL_n({\Q}_p).$  
\begin{itemize}
\item[(1)] Let $S_0$ be of type (M-3) or (M-5).  Then for any $k
  \ge n$ we have  
\begin{equation*}\begin{split}\beta_p(&H_{k+l+1},-B_l)\alpha_p(H_{k-l-1}
    \perp B_l,T)\\&=\sum_{i=0}^{l}(-1)^ip^{i(i-1)/2
    +i(n-2k+1)}C_{2l+1,i}\alpha_p(H_{k+l+1}, -\tilde B_{l,i} \perp T), 
\end{split}
\end{equation*}
where 
$C_{m,i}={\prod_{j=1}^i (p^{m+1-2j}-1) \over \prod_{j=1}^i (p^j-1)}$
for  an odd positive integer $m$ and an integer  $i$ such that $ i \le
(m-1)/2.$ 

\item[(2)] Let $S_0$ be of type (M-1),(M-2),(M-4), or (M-6). Put
$\epsilon=\epsilon(S_0)=-1$ or $1$ according as $S_0$ is of type (M-6) or not. 
Then for any $k \ge n$ we have 
\begin{equation*}\begin{split}\beta_p(&H_{k+l+1},-B_l)\alpha_p(H_{k-l-1}
    \perp B_l,T)\\
&=\sum_{i=0}^{l}(-1)^ip^{i(i-1)/2 +i(n-2k+1)}
C_{2l,i,\epsilon}
\alpha_p(H_{k+l+1}, -\tilde B_{l,i} \perp T),
\end{split}\end{equation*} 
where 
$C_{m,i,\epsilon}={(p^{m/2}-\epsilon)(p^{m/2-i}+\epsilon)\prod_{j=1}^{i-1} (p^{m-2j}-1) \over \prod_{j=1}^i (p^j-1)}$
for  an even positive integer $m$ and an integer  $i$ such that $ i \le
m/2,$ and $\epsilon=\pm 1.$

\item [(3)] $$\alpha_p(H_{k+l}, -H \perp
  T)=(1-p^{-(k+l)})(1+p^{-(k+l-1)})\alpha_p(H_{k+l-1}, T)$$  
\end{itemize}
\end{ksplemma}
\begin{proof}
By Proposition 2.2 of \cite{katsurada-amj99}, 
we have 
\begin{equation*}\begin{split}
\beta_p(H_{k+l+1},-B_l)&\alpha_p(H_{k-l-1} \perp B_l,T)\\
&=\sum_{i=0}^{2l+2}(-1)^ip^{i(i-1)/2 +i((n+2l+2)+1-(2k+2l+2))}\\
&\quad\quad\times \sum_{
 G \in GL_{2l+2}({\Z}_p) \backslash \pi_{2l+2,i}
 }\alpha_p(H_{k+l+1}, -B_l[G^{-1}] \perp T).
\end{split}\end{equation*} 
 We note that $\alpha_p(H_{k+l+1},-B_l[G^{-1}] \perp T)=0$ if $G \in
 \pi_{2l+2,i}$ with $i \ge l+1.$ Fix $i=0,1,...,l.$ Then  by
 Lemma 2.3 of \cite{IK-squared}, we have $-B_l[G^{-1}] \perp T \sim
 -\tilde B_{l,i} \perp T$ 
 if $G \in \pi_{l,i}$ and $B_l[G^{-1}] \in {\mathcal H}_{2l+2}({\Z}_p).$  
Furthermore, by Proposition 2.8 of  \cite{IK-squared} we have
 $$\#(GL_{2l+2}({\Z}_p) \backslash \{ G \in \pi_{2l+2,i} \ ; \
 B_l[G^{-1}] \in {\mathcal H}_{2l+2}({\Z}_p)\}) ={\prod_{j=1}^i
 (p^{2l+2-2j}-1) \over \prod_{j=1}^i (p^j-1)}.$$ 
This proves the assertion (1). Similarly, the assertion (2) can be
 proved. Now again by Proposition 
 2.2 of \cite{katsurada-amj99} we have 
$$\beta_p(H_{k+l},H)\alpha_p(H_{k+l-1} ,T)= \alpha_p(H_{k+l}, -H \perp T).$$
On the other hand, we have 
$$\beta_p(H_{k+l},H)=(1-p^{-(k+l)})(1+p^{-(k+l-1)})$$
(e.g. Lemma 9, \cite{Ki-dirichlet}.) Thus the assertion (3) holds. 
\end{proof}
Now for a non-degenerate half-integral matrix $B$ of degree $n$ over
${\Z}_p$ define a polynomial $\gamma_p(B;X)$ in $X$ by 
$$\gamma_p(B;X)=
\left\{
\begin{array}{ll}
(1-X)\prod_{i=1}^{n/2}(1-p^{2i}X^2)(1-p^{n/2}\xi_p(B)X)^{-1} & \ {\rm
  if} \ n \ {\rm is \ even} \\ 
(1-X)\prod_{i=1}^{(n-1)/2}(1-p^{2i}X^2) & \ {\rm if} \ n \ {\rm is \ odd}
\end{array}
\right.$$

\bigskip

For a half-integral matrix $B$ of degree over $\Z_p,$ let $(\bar W,\bar q)$ denote the quadratic space over $\Z_p/p\Z_p$ defined by the quadratic form $\bar q({\bf x})=B[{\bf x}] \ {\rm mod} \ p,$ and define the radical $R(\bar W)$ of $\bar W$ by
$$R(\bar W)=\{{\bf x} \in \bar W ; \bar B({\bf x},{\bf y})=0 \ {\rm for \ any } \ {\bf y} \in \bar W \},$$
where $\bar B$ denotes the associated symmetric bilinear form of $\bar q.$ 
We then put $l_p(B)= {\rm rank}_{\Z_p/p\Z_p} R(\bar W)^{\perp},$ where
$R(\bar W)^{\perp}$ is the orthogonal complement of $R(\bar W)$ in
$\bar W.$ Furthermore, in case $l_p(B)$ is even,  put $\bar
\xi_p(B)=1$ or $-1$ according as $R(\bar W)^{\perp}$ is hyperbolic or
not. Here we make the convention that $\xi_p(B)=1$ if $l_p(B)=0.$ We
note that $\bar \xi_p(B)$ is different from $\xi_p(B).$ 

\bigskip

\begin{ksplemma}\label{lemma2} 
\begin{itemize}
\item[(1)] 
Let $B$ be a half-integral matrix of degree $n$ over $\Z_p.$ Put $l=l_p(B).$  
Then we  have
\begin{equation*}
%\begin{split}
\beta_p(H_m%&
,B)=
%\\
%&\quad
(1-p^{-m}) 
%\left\{\begin{array}{ll}
%\begin{cases} 
(1+\bar \xi_p(B)p^{n-l/2-m})\prod_{j=0}^{n-l/2-1}(1-p^{2j-2m})
%& \ {\rm if} \ l \ {\rm is \ even} \\
%\end{split}
\end{equation*}
if $l$ is even,
\begin{equation*}
%\begin{split}
\beta_p(H_m,B)=\prod_{j=0}^{n-(l+1)/2}(1-p^{2j-2m}) 
%& \ {\rm if} \ l \ {\rm is \ odd.}
%\end{cases}.
%\end{split}
\end{equation*}
if $l $ is odd.

\item[(2)] Let $ T \in {\mathcal H}_n({\Z}_p) \cap GL_n({\Q}_p).$ 
Then
there exists a polynomial $F_p(T,X)$ such that  
$\alpha_p(H_m, T)=\gamma_p(T;p^{-m})F_p(T,p^{-m})$.
\end{itemize}
\end{ksplemma}
\begin{proof}
The assertion (1) follows from Lemma 9, \cite{Ki-dirichlet}. The
assertion (2) is well known (cf. \cite{Ki-dirichlet}).  
\end{proof}
\bigskip

Let $({} \ ,{} \ )_p$ be the Hilbert symbol over ${\Q}_p$ and $h_p$
the Hasse invariant (for the definition of the Hasse invariant, see
\cite{Ki}).  Let $B$ be a non-degenerate symmetric matrix of degree $n$
with entries in ${\Q}_p.$ We define
\begin{equation*}
\begin{cases} 
\eta_p(B)= h_p(B)(\det
B,(-1)^{(n-1)/2} \det B)_p & \quad\text{ if } n \text{ is odd}\\
\xi_p(B)=\chi_p((-1)^{n/2}\det B)&\quad \text{ if } n \text{ is even}.
\end{cases}
\end{equation*}
From now on we often
write $\xi(B)$ instead of $\xi_p(B)$ and so 
on if there is no fear of confusion.  For a non-degenerate
half-integral matrix $B$ of degree 
$n$ over ${\Z}_p$ put $D(B)=\det B$ and $d(B)={\rm ord}_p(D(B)).$ Further,
put $$\delta(B)= \left\{\begin{array}{ll} 2[(d(B)+1)/2] & \ {\rm if}
\ n \ {\rm is \ even} \\ d(B) & \ {\rm if} \ n \ {\rm is \ odd }
\end{array}\right..$$ 
Let $\nu(B)$ be the least integer $l$ such that $p^l B^{-1} \in
{\mathcal H}_n({\Z}_p).$   Further put $\xi'(B)=1+\xi(B)-\xi(B)^2$ for
a matrix $B$ of even degree.  
Then we have

\bigskip

\begin{kspproposition}\label{prop1}
Let $B_1=(b_1)$ and $B_2$ be non-degenerate  half-integral matrices of degree $1$ and $n-1,$ respectively over ${\Z}_p,$ and put $B=B_1 \perp B_2.$ Assume that $ {\rm ord}_p(b_1) \ge \nu(B_2)-1.$   
\begin{itemize}
\item[(1)]
Let $n$ be even. Then we have 
\begin{equation*}\begin{split}F_p(-(H_i &\perp pH_{l-i}) \perp
    B,p^{-(k+l)})\\
&= 
{1-\xi p^{n/2 -k} \over  1-p^{n-2k+1}}
F_p(-(H_i \perp pH_{l-i}) \perp B_2
,p^{-(k+l-1)})\\
&\quad+K(B)p^{l-i}F_p(-(H_i \perp pH_{l-i}) \perp B_2,p^{-(k+l)}),
\end{split}\end{equation*}
where $\xi=\xi(B),$ and $K(B)$ is a rational number depending only on $B.$

\item[(2)]
Let $n$ be odd. Then we have 
\begin{equation*}\begin{split}
F_p(-(H_i &\perp pH_{l-i}) \perp B,p^{-(k+l)})\\
&={1 \over  1-\tilde \xi p^{(n+1)/2-k}}F_p(-(H_i \perp pH_{l-i}) \perp
B_2,p^{-(k+l-1)})\\
&\quad + K(B)p^{l-i}F_p(-(H_i \perp pH_{l-i}) \perp B_2,p^{-(k+l)}),
\end{split}\end{equation*}
where $\tilde \xi=\xi(B_2),$ and $K(B)$ is a rational number depending
only on $B.$ Here we understand that $B_2$ is the empty matrix and
that we have $\xi=1$ if $n=1.$   
\end{itemize}
\end{kspproposition}
\begin{proof} (1) Let $n$ be even. We have ${\rm ord}_p(b_1) \ge \nu(-(H_i \perp pH_{l-i}) \perp B_2) -1.$ Thus by Theorem 4.1 of \cite{katsurada-amj99},
we have
\begin{equation*}\begin{split}
F_p(-(H_i &\perp pH_{l-i})  \perp B,p^{-(k+l)})\\ 
&= {1-\xi(l,i) p^{(n+2l)/2 -(k+l)} \over
  1-p^{n+2l+1-2(k+l)}}F_p(-(H_i \perp pH_{l-i}) \perp B_2,p^{-(k+l-1)})\\ 
&\quad +(-1)^{\xi(l,i)+1}\xi(l,i)'\tilde \eta(l,i)
{1-\xi(l,i)p^{(n+2l)/2+1-(k+l)} \over
  1-p^{n+2l+1-2(k+l)}}\\
&\quad\times(p^{(n+2l)/2-(k+l)})^
{\delta(l,i)-\tilde \delta(l,i) +\xi(l,i)^2}p^{\delta(l,i)/2}\\
&\quad\times F_p(-(H_i \perp pH_{l-i}) \perp B_2,p^{-(k+l)}),
\end{split}\end{equation*}
where $\xi(l,i)=\xi(-(H_i \perp pH_{l-i}) \perp B), \xi(l,i)'=\xi'(- (H_i \perp pH_{l-i}) \perp B),\tilde \eta(l,i)=\eta(-(H_i \perp pH_{l-i}) \perp B_2),
\delta(l,i)=\delta(-(H_i \perp pH_{l-i}) \perp B),$ and $\tilde
\delta(l,i)=\delta(-(H_i \perp pH_{l-i}) \perp B_2).$ 
We note that $\xi(l,i),\xi(l,i)'$ and $\tilde \eta(l,i)$ are
independent of $l$ and $i,$ and they are equal to $\xi, \xi',$ and
$\eta(B_2),$ respectively. 
Furthermore, we have $\delta(l,i)=2l-2i+2[({\rm ord}_p( \det
T)+1)/2]$ and $\tilde \delta(l,i)=2l-2i+{\rm ord}_p(\det \hat T).$   
Thus the assertion holds. Similarly, the assertion holds in case $n$ is odd.

\end{proof}

\bigskip

\begin{kspproposition}\label{prop2}
Let $S_0$ and the others be as in Lemma \ref{lemma1}. $T=b_1 \perp b_2 \perp ... \perp b_n$ with ${\rm ord}_p(b_1) \ge
    {\rm ord}_p(b_2) \ge ...\ge {\rm ord}_p(b_n).$ Put $\hat T= b_2
    \perp ...\perp b_n.$   
 \begin{itemize}
\item[(1)] Assume that  $n+ \deg S_0$ is even. Put $K(S_0,T)={1-p^{n-2k}
    \over  1-p^{n/2-k} \xi}K(-S_0 \perp T),$ where $\xi=\xi(-S_0 \perp
  T),$ and $K(-S_0 \perp T)$ is the rational number in Proposition
  \ref{prop1} .  
Then we have 
\begin{equation*}\begin{split}\alpha_p(H_{k+l+1},& -\tilde B_{l,i}
    \perp T)\\
&={(1-p^{-(k+l+1)})(1+p^{-(k+l)}) \over 1-p^{n-2k+1}}
\alpha_p(H_{k+l},-\tilde B_{l,i} \perp \hat T)\\
&\quad +p^{l-i}K(S_0,T)
\alpha_p(H_{k+l+1}, -\tilde B_{l,i} \perp \hat T).
\end{split}\end{equation*}

\item[(2)] Assume that $n+ \deg S_0$ is odd. Put 
$K(S_0,T)=(1-p^{(n-1)/2-k} \tilde \xi)K(-S_0 \perp T),$ where $\tilde
\xi=\xi(-S_0 \perp \hat T),$ and $K(-S_0 \perp T)$ is the rational number
in Proposition \ref{prop1}.  Then we have 
\begin{equation*}\begin{split}\alpha_p(H_{k+l+1},& -\tilde B_{l,i} \perp T)\\
&={(1-p^{-(k+l+1)})(1+p^{-(k+l)}) \over 1-p^{n-2k+1}}
\alpha_p(H_{k+l},-\tilde B_{l,i} \perp \hat T)\\
&\quad+p^{l-i}K(S_0,T)
\alpha_p(H_{k+l+1}, -\tilde B_{l,i} \perp \hat T).
\end{split}\end{equation*}

\end{itemize}

\end{kspproposition}

\begin{proof}
By (1) of Proposition \ref{prop1} and (2) of Lemma \ref{lemma2}, we
have 
\begin{eqnarray*}
\alpha_p(H_{k+l+1},-\tilde B_{l,i} \perp T)&=&
\gamma_p(-\tilde B_{l,i}
\perp T, p^{-(k+l+1)})F_p(-\tilde B_{l,i} \perp T,p^{-(k+l+1)})\\ 
&=&\gamma_p(-\tilde B_{l,i} \perp T, p^{-(k+l+1)})\\
&&\quad\times [{1-\xi p^{n/2 -k} \over
  1-p^{n-2k+1}}F_p(-\tilde B_{l,i} \perp \hat
T,p^{-(k+l)})\\ 
&&+p^{l-i}K(-S_0 \perp T)F_p(-\tilde B_{l,i} \perp \hat
T,p^{-(k+l+1)})].
\end{eqnarray*}
We note that 
\begin{equation*}\begin{split}
\gamma_p(-\tilde B_{l,i} \perp \hat T,& p^{-(k+l)})\\
&={1-p^{n/2-k} \xi
  \over  (1-p^{-(k+l+1)})(1+p^{-(k+l)})} \gamma_p(-\tilde B_{l,i} \perp
T, p^{-(k+l+1)}),
\end{split}\end{equation*} 
and 
$$\gamma_p(-\tilde B_{l,i} \perp \hat T, p^{-(k+l+1)})={1-p^{n/2-k} \xi
  \over  1-p^{n-2k}} \gamma_p(-\tilde B_{l,i} \perp T, p^{-(k+l+1)}).$$ 
Thus the assertion (1) holds. 

Now by (2) of Proposition \ref{prop1} and (2) of Lemma \ref{lemma2}, we have 
\begin{eqnarray*}\alpha_p(H_{k+l+1},-\tilde B_{l,i} \perp T)
&=&\gamma_p(-\tilde B_{l,i}
\perp T, p^{-(k+l+1)})F_p(-\tilde B_{l,i} \perp T,p^{-(k+l+1)})\\
&=&\gamma_p(-\tilde B_{l,i} \perp T, p^{-(k+l+1)})\\
&&\quad\times [{ 1  \over 1- \tilde \xi p^{(n+1)/2-k}}
F_p(-\tilde B_{l,i} \perp \hat T,p^{-(k+l)})\\ 
&&+p^{l-i}K(-S_0 \perp T)   F_p(-\tilde B_{l,i} \perp \hat
T,p^{-(k+l+1)})].
\end{eqnarray*}
We note that 
\begin{equation*}\begin{split}
\gamma_p(-\tilde B_{l,i} \perp \hat T, p^{-(k+l)})&=
{1-p^{n+1-2k}
  \over  (1-p^{-(k+l+1)})(1+p^{-(k+l)})(1- \tilde \xi p^{(n+1)/2-k}) }\\
&\quad \times
\gamma_p(-\tilde B_{l,i} \perp T, p^{-(k+l+1)}),
\end{split}\end{equation*}
and 
$$\gamma_p(-\tilde B_{l,i} \perp \hat T, p^{-(k+l+1)})={1 \over 1-
  p^{(n-1)/2-k} \tilde \xi} \gamma_p(-\tilde B_{l,i} \perp T,
p^{-(k+l+1)}).$$ 
Thus the assertion (2) holds.

\end{proof}

\begin{ksprem}
\begin{rm}
In  the above theorem, $K(S_0,T)$ can be expressed explicitly in terms
of the invariants of $T.$ 
\end{rm}
\end{ksprem}

\begin{kspproposition}\label{prop3}{
Let $S_0,T$ and $\hat T$ and the others be  as in Proposition \ref{prop2}. 
 \begin{itemize}
\item[(1)] Assume that  $S_0$ is of type (M-3) or (M-5). Then for any
 non-negative integer $l \le k-1$ we have   
\begin{equation*}\begin{split}
\alpha_p(H_{k-l-1} \perp B_l,T)&= (1-p^{n-2k+1})^{-1}\\
&\quad\times\{(1-p^{-2k+2l+2}) \alpha_p(H_{k-l-2} \perp B_l,\hat T)\\
&\quad\quad +
p^{n-2k+1}(p^{2l}-1) \alpha_p(H_{k-l-1} \perp B_{l-1},\hat T) \}\\
&\quad+ p^l K(S_0,T) \alpha_p(H_{k-l-1} \perp B_l,\hat T),
\end{split}\end{equation*}
where $K(S_0,T)$ is the rational number in Proposition \ref{prop1}.
In particular, if $n=1$, for  a non-zero element $T$ of ${\Z}_p$, we have
 $$\alpha_p(H_{k-l-1} \perp B_l, T)= 1 +  c  p^l,$$
 where $c=c(S_0,T)$ is the rational number determined by $T$ and $S_0.$ 

\item[(2)] Assume that $S_0$ is of type (M-1),(M-2),(M-4) or (M-6). Put
  $l'=l+1$ or $l$ according as $S_0$ is of type (M-6) or not. 
Let $\epsilon=\epsilon(S_0)$ be as in Lemma \ref{lemma1}, and $\bar \xi
=\bar \xi(S_0).$ Put $\epsilon= 
-1$ or $1$ according as $S_0$ is of type (M-6) or not. Then  for
non-negative integer $l 
\le k-1$ we have   
\begin{equation*}\begin{split}
\alpha_p(&H_{k-l-1} \perp B_l ,T)=(1-p^{n-2k+1})^{-1}\\ 
&\quad\times\{(1-p^{-k+l'+1} \bar \xi) (1+p^{-k+l'} \bar \xi)\alpha_p(H_{k-l-2}
\perp B_l,\hat T)\\
&\quad\quad + 
p^{n-2k+1}(p^{l'}-\epsilon)(p^{l'-1}+\epsilon)\alpha_p(H_{k-l-1} \perp
B_{l-1},\hat T)\}\\
&\quad+K(S_0,T)  p^l \alpha_p(H_{k-l-1} \perp B_l,\hat T),
\end{split}\end{equation*}
where $K(S_0,T)$ is 
the  rational number in Proposition \ref{prop1}.
%\newpage
In particular, if $n=1$, for  a non-zero element $T$ of ${\Z}_p$, we have
 $$\alpha_p(H_{k-l-1} \perp B_l, T)= 1 +  c p^l,$$
 where $c=c(S_0,T)$ is a rational number determined by $T$ and $S_0.$ 
Throughout (1) and (2), we understand   $\alpha_p(H_{k-l-2} \perp B_l,\hat T)=1$ if $l=k-1.$ 
\end{itemize}
}

\end{kspproposition}

\begin{proof}(1) 
First let $n+ \deg S_0$ be even.  Then by (1) of Proposition \ref{prop2}
and (1) of Lemma \ref{lemma1}, we have  
\begin{equation*}\begin{split}
\beta_p(&H_{k+l+1},-B_l)\alpha_p(H_{k-l-1} \perp B_l, T)\\
&=\sum_{i=0}^l (-1)^ip^{i(i-1)/2+i(n-2k+1)}C_{2l+1,i}\\
&\quad\times
\{{(1-p^{-(k+l+1)})(1+p^{-(k+l)})
  \over 1-p^{n-2k+1} } \alpha_p(H_{k+l},-\tilde B_{l,i} \perp \hat T)\\
&\quad\quad+\alpha_p(H_{k+l+1},-\tilde B_{l,i} \perp \hat T)p^{l-i}K(S_0 ,T)\}\\ 
&={(1-p^{-(k+l+1)})(1+p^{-(k+l)}) \over 1-p^{n-2k+1} }
\\
&\quad\times[\sum_{i=0}^l
(-1)^ip^{i(i-1)/2+i(n-2k+2)}
C_{2l+1,i}\alpha_p(H_{k+l},-\tilde B_{l,i} \perp \hat T)\\
&\quad\quad+\sum_{i=0}^l
(-1)^ip^{i(i-1)/2+i(n-2k+2)}(p^{-i}-1)C_{2l+1,i}\alpha_p(H_{k+l},-\tilde
B_{l,i} \perp \hat T)]\\ 
&\quad+\sum_{i=0}^l
(-1)^ip^{i(i-1)/2+i(n-2k+1)}C_{2l+1,i}\alpha_p(H_{k+l+1},-\tilde
B_{l,i} \perp \hat T)p^{l}K(S_0,T).
\end{split}\end{equation*} 
By (1) of Lemma \ref{lemma1} and (1) of Lemma \ref{lemma2}, we have
\begin{equation*}\begin{split}
(1-&p^{-(k+l+1)})(1+p^{-(k+l)})\\&\quad\times\sum_{i=0}^l
(-1)^ip^{i(i-1)/2+i(n-2k+2)}C_{2l+1,i}\alpha_p(H_{k+l},-\tilde B_{l,i}
\perp \hat T)\\
&=(1-p^{-(k+l+1)})(1+p^{-(k+l)})\beta_p(H_{k+l},-B_l)\alpha_p(H_{k-l-2}
\perp B_l, \hat T) \\ 
&=(1-p^{2l+2-2k})\beta_p(H_{k+l+1},-B_l)\alpha_p(H_{k-l-2} \perp B_l,
\hat T) 
\end{split}\end{equation*}
and  
\begin{equation*}\begin{split}
\sum_{i=0}^l
(-1)^ip^{i(i-1)/2+i(n-2k)}C_{2l+1,i}&\alpha_p(H_{k+l+1},-\tilde
B_{l,i} \perp \hat T)\\
&=\beta_p(H_{k+l+1},-B_l)\alpha_p(H_{k-l-1} \perp B_l, \hat T).
\end{split}\end{equation*}
Furthermore, again by (1) and (3) of Lemma \ref{lemma1}, and (1) of
Lemma \ref{lemma2},  we have   
%\vfill
%\newpage
\begin{equation*}\begin{split}
(1-&p^{-(k+l+1)})(1+p^{-(k+l)})\sum_{i=0}^l
(-1)^ip^{i(i-1)/2+i(n-2k+2)} (p^{-i}-1)\\
&\quad\quad\quad\quad\quad\quad\quad\quad\quad\quad
\times C_{2l+1,i}\alpha_p(H_{k+l},-\tilde B_{l,i}
\perp \hat T)\\ 
&=(1-p^{-(k+l+1)})(1+p^{-(k+l)})p^{n-2k+1}(p^{2l}-1)\\
&\quad\quad\times \sum_{i=1}^l (-1)^{i-1}p^{(i-2)(i-1)/2+(i-1)(n-2k+2)}
C_{2l-1,i-1}\\
&\quad\quad\quad\quad\quad
\times\alpha_p(H_{k+l},-\tilde B_{l-1,i-1} \perp -H \perp \hat
T)\\
&=(1-p^{-(k+l+1)})(1+p^{-(k+l)})p^{n-2k+1}(p^{2l}-1)\\
&\quad\quad\times \sum_{i=1}^l (-1)^{i-1}p^{(i-2)(i-1)/2+(i-1)(n-2k+2)}\\
&\quad\quad\times
C_{2l-1,i-1}(1-p^{-(k+l)})(1+p^{-(k+l-1)})\\
&\quad\quad\times\alpha_p(H_{k+l-1},-\tilde
B_{l-1,i-1} \perp \hat T)\\
&=p^{n-2k+1}(p^{2l}-1)(1-p^{-(k+l+1)})(1-p^{-2(k+l)})(1+p^{-(k+l-1)})\\
&\quad\quad\quad\quad\quad\quad
\times \beta_p(H_{k+l-1},-B_{l-1})\alpha_p(H_{k-l-1} \perp
B_{l-1}, \hat T)\\
&=p^{n-2k+1}(p^{2l}-1)\beta_p(H_{k+l+1},-B_{l-1})\alpha_p(H_{k-l-1}
\perp B_{l-1}, \hat T).
\end{split}\end{equation*}
This proves  the assertion (1) in case $n +\deg S_0$ is odd. 
Next again by
(2) of Proposition \ref{prop2} and (1) of Lemma \ref{lemma1}, the
assertion (1) can be proved  in case $n+\deg S_0$ is odd. 
%\vfill
%\newpage

(2)  First let $n+ \deg S_0$ be even.  Then by  (1) of Proposition
\ref{prop2} and (2) of Lemma \ref{lemma1}, we have  
\begin{equation*}\begin{split}
\beta_p(&H_{k+l+1},-B_l)\alpha_p(H_{k-l-1} \perp B_l, T)\\
&=\sum_{i=0}^{l'}
(-1)^ip^{i(i-1)/2+i(n-2k+1)}C_{2l',i,\epsilon}\\
&\quad\times\{{(1-p^{-(k+l'+1)})(1+p^{-(k+l')})  \over 1-p^{n-2k+1} } 
\alpha_p(H_{k+l},-\tilde B_{l,i} \perp \hat T)\\
&\quad\quad\quad+\alpha_p(H_{k+l+1},-\tilde B_{l,i} \perp \hat
T)p^{l-i}K(S_0,T)\}.
\end{split}\end{equation*}
We evaluate this further as

\begin{multline*} 
{(1-p^{-(k+l'+1)})(1+p^{-(k+l')}) \over 1-p^{n-2k+1} }\\
\shoveleft{\quad \quad\times[\sum_{i=0}^{l'}
(-1)^ip^{i(i-1)/2+i(n-2k+2)}
C_{2l',i,\epsilon}\alpha_p(H_{k+l},-\tilde B_{l,i} \perp \hat T)}\\
\shoveleft{\quad\quad +\sum_{i=0}^{l'}
(-1)^ip^{i(i-1)/2+i(n-2k+2)}(p^{-i}-1)C_{2l',i,\epsilon}\alpha_p(H_{k+l},-\tilde
B_{l,i} \perp \hat T)]}\\  
\shoveleft{\quad+\sum_{i=0}^{l'}
(-1)^ip^{i(i-1)/2+i(n-2k)}C_{2l',i,\epsilon}
\alpha_p(H_{k+l+1},-\tilde
B_{l,i} \perp \hat T)p^{l}K(S_0,T).}
\end{multline*}
By (1) and (3) of Lemma \ref{lemma1} and (1) of Lemma \ref{lemma2}, we have
\begin{equation*}\begin{split}
(1&-p^{-(k+l'+1)})(1+p^{-(k+l')})\\
&\quad\times \sum_{i=0}^l
(-1)^ip^{i(i-1)/2+i(n-2k+2)}C_{2l',i,\epsilon}\alpha_p(H_{k+l},-\tilde B_{l,i}
\perp \hat T)\\
&=(1-p^{-(k+l'+1)})(1+p^{-(k+l')})\beta_p(H_{k+l},-B_l)\alpha_p(H_{k-l-2}
\perp B_l, \hat T) \\
&=(1-\bar \xi p^{l'+1-k})(1+\bar \xi
p^{l'-k})\beta_p(H_{k+l+1},-B_l)\alpha_p(H_{k-l-2} \perp B_l, \hat T) 
\end{split}\end{equation*}
and 
\begin{equation*}\begin{split}
\sum_{i=0}^l
(-1)^ip^{i(i-1)/2+i(n-2k)}&C_{2l',i,\epsilon}\alpha_p(H_{k+l+1},-\tilde
B_{l,i} \perp \hat T)\\ 
&=\beta_p(H_{k+l+1},-B_l)\alpha_p(H_{k-l-1} \perp B_l, \hat T).
\end{split}\end{equation*}
Furthermore, again by (1) of Lemma \ref{lemma1}, and (1) of Lemma
\ref{lemma2},  we have   
\begin{equation*}\begin{split}
(1-&p^{-(k+l'+1)})(1+p^{-(k+l')})\sum_{i=0}^{l'}
(-1)^ip^{i(i-1)/2+i(n-2k+2)} (p^{-i}-1)\\ 
&\quad\quad\times C_{2l',i,\epsilon}\alpha_p(H_{k+l},-\tilde B_{l,i}
\perp \hat T)\\
&=(1-p^{-(k+l'+1)})(1+p^{-(k+l')})p^{n-2k+1}(p^{l'}-\epsilon)(p^{l'-1}+\epsilon)\\
&\quad\times \sum_{i=1}^{l'} (-1)^{i-1}p^{(i-2)(i-1)/2+(i-1)(n-2k+2)}\\
&\quad\quad\times C_{2l'-2,i-1,\epsilon}\alpha_p(H_{k+l},-\tilde
B_{l-1,i-1} \perp -H \perp \hat  T),
\end{split}\end{equation*}
which can be transformed into
\begin{equation*}\begin{split}
(1-&p^{-(k+l'+1)})(1+p^{-(k+l')})p^{n-2k+1}(p^{l'}-\epsilon)(p^{l'-1}+\epsilon)\\
&\quad\times \sum_{i=1}^l
(-1)^{i-1}p^{(i-2)(i-1)/2+(i-1)(n-2k+2)}C_{2l'-2,i-1,\epsilon} \\
&\quad\quad\times
(1-p^{-k+l'})(1+p^{-(k+l'-1)})\alpha_p(H_{k+l-1},-\tilde
B_{l-1,i-1} \perp \hat T)\\ 
&=p^{n-2k+1}(p^{l'}-\epsilon)(p^{l'-1}+\epsilon)(1-p^{-(k+l'+1)})(1-p^{-2(k+l')})\\
&\quad\times (1+p^{-(k+l'-1)})\beta_p(H_{k+l-1},-B_{l-1})\alpha_p(H_{k-l-1} \perp
B_{l-1}, \hat T)\\
&=p^{n-2k+1}(p^{l'}-\epsilon)(p^{l'-1}+\epsilon)\\
&\quad\quad\times\beta_p(H_{k+l+1},-B_{l-1})\alpha_p(H_{k-l-1}  
\perp B_{l-1}, \hat T).
\end{split}\end{equation*} 
This proves  the assertion (2) in case $n +\deg S_0$ is odd. Next again
by (2) of Proposition \ref{prop2} and Lemma \ref{lemma1}, the
assertion (2) can be proved  in case $n+\deg S_0$ is odd. 
\end{proof}
\begin{proof}[Proof of Theorem \ref{densitypolynomial}]
We prove the assertion by induction on $n.$ 
The assertion for $n=1$ follows from (2) of
Proposition\ref{prop3}. Let $n \ge 2$ and assume that the assertion
holds for $n-1.$  Then 
by the induction assumption  we have 
$$\alpha_p(H_{s-t-1} \perp B_t, \hat T)=\sum_{j=0}^{n-1} a_j p^{tj},$$
and
$$\alpha_p(H_{s-t-2} \perp B_t, \hat T)=\sum_{j=0}^{n-1} a_j' p^{tj},$$
where $a_j=a_j(s,S_0,\hat T)$ and $a_j'(s-1,S_0,\hat T)$ in Theorem \ref{densitypolynomial}. We may assume that $T= b_1 \perp b_2 \perp ... \perp b_n$ with 
${\rm ord}_p(b_1) \ge {\rm ord}_p(b_2) \ge ... \ge {\rm ord}_p(b_n).$
First assume that $S_0$ is of type (M-3) or (M-5). Thus by  Proposition \ref{prop3} we have 
\begin{equation*}\begin{split}
\alpha_p(H_{k-l-1} \perp B_l, T)&={1-p^{-2k+2l+2} \over 1-p^{n-2k+1}}
\alpha_p(H_{k-l-2} \perp B_l,\hat T)\\
&\quad + {p^{n-2k+1}(p^{2l}-1) \over 1-p^{n-2k+1}}\alpha_p(H_{k-l-1} \perp
B_{l-1},\hat T)\\
&\quad + p^l K(S_0,T)   \alpha_p(H_{k-l-1} \perp B_l,\hat T)
\end{split}\end{equation*}
which is equal to
 
\begin{equation*}\begin{split}{1-p^{-2k+2l+2} \over 1-p^{n-2k+1}}
    \sum_{j=0}^{n-1} a_{j}' p^{lj} &+ {p^{n-2k+1}(p^{2l}-1) \over
      1-p^{n-2k+1}} \sum_{j=0}^{n-1} 
 a_j'p^{(l-1)j}\\
&+p^l  K(S_0,T) \sum_{j=0}^{n-1} a_{j}p^{lj}.
\end{split}\end{equation*}
For $ 0 \le j \le n-1$ put
$$M(j)= {1-p^{-2k+2l+2} \over 1-p^{n-2k+1}}  a_j' p^{lj}  +
{p^{n-2k+1}(p^{2l}-1) \over 1-p^{n-2k+1}}  a_j' p^{(l-1)j}+p^l  K(S_0,T)
a_j p^{lj}.$$ 
Then for $j \le n-2, M(j)$ is a polynomial in $p^l$ of degree at most $n-1.$ 
On the other hand,  
\begin{eqnarray*}M(n-1)&=&\frac{1-p^{-2k+2l+2}} {1-p^{n-2k+1}}
  a_{n-1}' p^{l(n-1)}+  \frac{p^{n-2k+1}(p^{2l}-1)} {1-p^{n-2k+1}}
  a_{n-1}' p^{(l-1)(n-1)}\\
&&+ a_{n-1}p^l  K(S_0,T) p^{l(n-1)}\\
&=&a_{n-1}'{1-p^{-2k+2} \over 1-p^{n-2k+1}}
p^{l(n-1)}+a_{n-1}K(T)p^{ln}.
\end{eqnarray*} 
Thus $\alpha_p(H_{k-l-1} \perp B_l, T)$ is a polynomial in $p^l$ of
degree at most $n.$  This proves the assertion in case (M-3) or
(M-5). Similarly, the assertion can be proved in the remaining case.   

\bigskip

\begin{ksprem}
\begin{rm}

A more careful analysis shows that we have  $a_0(k,S_0,T)=1$ in the
above theorem. 

\end{rm}
\end{ksprem}

{\bf Corollary to Theorem \ref{densitypolynomial}} {\it Let the
  notation be as above. For 
  any $n$-tuple $(l_1,l_2,....,l_n)$ of complex numbers, put
  $\mu(l_1,...,l_n)= \prod_{1 \le j \le i \le n} (p^{l_i}-p^{l_j}).$
  Then for any integers $0 \le l_1< ...<l_{n+2} \le k$ and $T \in
  {\mathcal H}_n({\Z}_p) \cap GL_n({\Q}_p)$ we have 
$$\sum_{j=1}^{n+2}
(-1)^{j-1}\mu(l_1,...,l_{j-1},l_{j+1},...,l_{n+2})\alpha_p(H_{k-l_j
  -1} \perp B_{l_j}, T)=0.$$} 
\end{proof}
\begin{ksptheorem}\label{theorem2}{Let $k \ge n+1.$ Let $n+1$ integers 
$0 \le l_1 ...<l_{n+2} \le k$ be given, let
  $\lambda_1,...,\lambda_{n+1}$ be  rational numbers such that  
$$\sum_{j=1}^{n+1} \lambda_j \alpha_p(H_{n-l_j+1} \perp B_{l_j+k-n-2},T)=0$$
for any $T \in {\mathcal H}_n({\Z}_p) \cap GL_n({\Q}_p)$. Then we have
$\lambda_1=....=\lambda_{n+1}=0.$ } 
\end{ksptheorem}
\begin{proof} We prove the assertion by induction on $n.$ The assertion
 for $n=1$ follows from   Proposition \ref{prop3}.  Let $n \ge 2,$ and
 assume that the assertion holds for $n-1.$ The above relation holds
 for $T= p^{2r} \perp \hat T$ with any  
 integer $r$ and $\hat T  \in {\mathcal H}_{n-1}({\Z}_p) \cap GL_{n-1}({\Q}_p).$ 
Then by Proposition \ref{prop3}, 
%\begin{equation*}\begin{split}
\begin{multline*}
\sum_{l=1}^{n+1} \lambda_l\{(1-p^{2l-2n-2})\alpha_p(H_{n-l} \perp
B_{l+k-n-2}, \hat T)\\
\shoveleft{\quad\quad +p^{n-2k+1}(p^{2l+2k-2n-4}-1)\alpha_p(H_{n-l+1} \perp
B_{l+k-n-3},\hat T)\}}\\ 
\shoveleft{
\quad\quad +p^{(n-2k+1)r} w(\hat T)\sum_{l=1}^{n+1} \lambda_l
p^{l+k-n-2}\alpha_p(H_{n-l+1} \perp B_{l+k-n-2}, \hat T)=0,}
\end{multline*}
%\end{split}\end{equation*}
where $w(\hat T)$ is a certain rational number depending only on $\hat T.$ 
Thus by taking the limit $r \rightarrow \infty$ we obtain 
\begin{multline*}
\sum_{l=1}^{n+1} \lambda_l\{(1-p^{2l-2n-2})\alpha_p(H_{n-l} \perp
B_{l+k-n-2}, \hat T)\\
%\shoveleft
{\quad\quad+p^{n-2k+1}(p^{2l+2k-2n-4}-1) \}\alpha_p(H_{n-l+1} \perp
B_{l+k-n-3},\hat T)=0 \qquad (*)}
\end{multline*} and  
$$\sum_{l=1}^{n+1} \lambda_l p^{l+k-n-2}\alpha_p(H_{n-l+1} \perp
B_{l+k-n-2}, \hat T)=0 \qquad (**).$$ 
Rewriting (*) we have 
\begin{equation*}\begin{split}\sum_{l=1}^{n} (\lambda_l(1-p^{2l-2n-2})
&+\lambda_{l+1}p^{n-2k+1}(p^{2l+2k-2n-2}-1))\\
&\quad\times \alpha_p(H_{n-l} \perp
B_{l+k-n-2},\hat T)\}=0 .
\end{split}\end{equation*}
Thus by the induction hypothesis, we have 
$$\lambda_l(1-p^{2l-2n-2}) +\lambda_{l+1}p^{n-2k+1}(p^{2l+2k-2n-2}-1)=0$$ 
for any $l=0,1,...,n.$ 
In particular 
$$\lambda_n= - { p^{n-2k+3}(p^{2k-2}-1) \over p^2-1}\lambda_{n+1} \qquad (***).$$
On the other hand, by the Corollary to Theorem \ref{densitypolynomial} we have 
$$\sum_{l=1}^{n+1} (-1)^{l-1} \mu_l \alpha_p(H_{n-l+1} \perp
B_{l+k-n-2}, \hat T)=0 \qquad (****),$$ 
where $\mu_l= \mu(k-n-1,...,l+k-n-3,l+k-n-1,...,k-1).$
By (**) and (****), and the induction hypothesis, we have
$$\lambda_l=(-1)^{l-n-1}{\mu_l \over \mu_{n+1}}p^{-l+n+1}
\lambda_{n+1},$$ 
and in particular
$$\lambda_n=-{\mu_n \over \mu_{n+1}}p \lambda_{n+1}=-{p^n -1 \over
  p-1} \lambda_{n+1}.  
\qquad (*****)$$
If $\lambda_{n+1} \not=0,$ (*****) contradicts (***), since $n \ge 2$
and $k \ge n+1.$ Thus we have $\lambda_{n+1}=0$ and therefore
$\lambda_l=0$ for any $l=1,...,n+1.$ This completes the induction.  
\end{proof}
We can now prove Theorem \ref{maintheorem}:
We notice first that the genera of lattices of level $p$ on the space
of the given lattice are represented by lattices $L^{(i)}$ whose
$p$-adic completions have a Gram matrix that is $\Z_p$-equivalent to 
$H_{k-i-1}\perp p H_i \perp S_0$ with  a fixed $S_0$ of degree $2$ as
in Theorem \ref{densitypolynomial}. Altogether there are $k \ge n+1$
such genera.

As a consequence of Siegel's theorem one sees that the linear
independence of any $n+1$ of the degree $n$ theta series of the genera
of the $L^{(i)}$ is implied by the linear independence of the
corresponding $p$-adic local density functions $T\mapsto
\alpha_p(L^{(i)},T)$ stated in Theorem \ref{theorem2} (notice that the
restriction $k\ge n+1$ implies that for primes $\ell\ne p$ the $\ell$-adic
completion of the lattices $L^{(i)}$ splits off an orthogonal sum of at
least $n$ unimodular hyperbolic planes so that every even $T \in
M_n^{\rm sym}(\Z)$ is represented at all $\ell$-adic completions and
the product of the $\alpha_\ell(L^{(i)},T)$ is nonzero).

Since all the genus theta series are (by Siegel's theorem) in the space
of Eisenstein series associated to zerodimensional boundary components
(cusps) and since there are $n+1$ such cusps in the case of prime
level,
it is clear that both types of genus theta series generate the full
space of Eisenstein series.
  
\medskip
\begin{kspcorollary}\label{hecke_expression}

Let $L$ be a lattice on the quadratic space $V$ over $q$ of level
$p$ as in Theorem \ref{maintheorem} and put 
$F=\vartheta^{(n)}(\text{gen}(L))$. Then the modular form $F\vert_kT(\ell)$
can be expressed as a linear combination of theta series of positive
definite lattices of level  $p$ on $V$ for all primes $\ell \ne p.$
\end{kspcorollary}
\begin{proof}
This is clear from Theorem \ref{theorem2} and Theorem \ref{heckeaction_indefinite}.
\end{proof}

\begin{ksprem}
\begin{rm}
\begin{itemize}
\item[a)] The result of Theorem \ref{maintheorem} is more generally true in the
case of square free level $N$, in which case the dimension of the
space spanned by the genus theta series becomes $(n+1)^{\omega(N)}$
where $\omega(N)$ is the number of primes dividing $N$; one has then a
basis of genus theta series if one considers  $(n+1)^{\omega(N)}$
genera of lattices on the 
same quadratic space $V$ such that for each $p$ dividing $n$ one has
$n+1$ local integral equivalence classes. In that case our proof given above
requires the restriction that the anisotropic kernel of the quadratic space
under consideration has dimension at most $2$. Moreover we can not
guarantee the holomorphy of the indefinite genus theta series if the
character is trivial (i.e., if the underlying quadratic space has
square discriminant). One proceeds in the proof as above, adding an
induction on the number of primes $\omega(N)$ dividing $N$.
\item[b)] A different (and much shorter) proof of Theorems
  \ref{densitypolynomial} and \ref{theorem2} has been communicated to
  us by Y. Hironaka and F. Sato \cite{hir_sa}. The proof given here
  gives a little more information (e.g. explicit recursion relations)
  than theirs. The proof of Hironaka and Sato removes the restriction
  on the anisotropic kernel mentioned above (if one strengthens the
  condition on $n$ to $n+1<k$ in the new cases) and provides also a
  version for levels that are not square free. The application of that
  version to the study of the space of Eisenstein series generated by
  the genus theta series in the case of arbitrary level will be the
  subject of future work.
\end{itemize}
\end{rm}
\end{ksprem}
 \section{Connection with Kudla's matching principle}
In Section 4 we have seen that the Hecke operator $T(p)$ can
provide a connection between theta series for lattices in positive 
definite quadratic space $(V_1,q_1)$ and in a related indefinite
quadratic space $(V_2,q_2)$. Such a connection has recently been
observed in a different setup by Kudla \cite{kudla-compositio03}. 
We sketch his approach
briefly in order to study the relation to our construction, for details
we refer to \cite{kudla-compositio03}, Section 4.1.
 \vspace{0.3cm}\\
Let $(V_1,q_1)$ be a positive definite quadratic space over ${\mathbb Q}$ of
dimension $m$ and discriminant $d$, let $(V_2,q_2)$ be a space of 
the same dimension $m$ and discriminant $d$, but of signature 
$(m-2,2)$. We fix $n > 0$ and an additive character $\psi$ of ${\mathbb Q}_
{\mathbb A}$. Consider the oscillator representations $\omega_1= \omega_{1,\psi}$
of $\widetilde{\rm Sp}_n({\mathbb A}) \times O_{(V_1,q_1)}({\mathbb A})$ on
the Schwartz space $S((V_1({\mathbb A}))^n)$ and $\tilde{\omega}$ of
$\widetilde{\rm Sp}_n({\mathbb A}) \times O_{(V_2,q_2)}({\mathbb A})$ on
$S((V_2({\mathbb A}))^n)$, where $\widetilde{\rm Sp}_n({\mathbb A})$
denotes the usual metaplectic double cover of the adelic symplectic group
${\rm Sp}_n({\mathbb A})$.
 \vspace{0.3cm}\\
For $j = 1,2$ we have then for $\varphi \in S((V_j({\mathbb A}))^n)$ the 
theta kernel
 $$\begin{array}{c}
 \theta(\tilde{g},h_j;\varphi_j) = \Sum_{x \in V_j({\mathbb Q})}
  \omega_j(\tilde{g}) \varphi_j(h_j^{-1} x)\vspace{0.2cm}\\
 (\tilde{g} \in \widetilde{\rm Sp}_n({\mathbb A}),\: h_j \in O_{(V_j,q)} ({\mathbb A})).
 \end{array}$$
and the theta integral
 $$I(\tilde{g}; \varphi_j) = \int_{O_{(V_j,q_j)}({\mathbb Q})\setminus
  O_{(V_j,q_j)}({\mathbb A})} 
  \theta(\tilde{g},{h}_j,\varphi_j) dh_j $$
which (under our conditions) is absolutely convergent for $j = 1$ and for
$j = 2$ if $V_2$ is anisotropic or $m > n+2$.
 \vspace{0.3cm}\\
Let now $L_j$ be a lattice on $V_j$ and assume $\varphi_j$ to be factored as
$\varphi_j = \prod_{v} \varphi_{j,v}$ over all places $v$ of ${\mathbb
  Q}$, where 
 $\varphi_{j,p} = {\mathbf 1}_{L_{j,p}}$ is the characteristic function 
of the lattice
$L_{j,p}$ in the ${\mathbb Q}_p$-space $V_{j,p}$ for all finite primes $p$.
Then for $\varphi_{1,\infty}({\mathbf x}) = \exp(-2\pi~{\rm tr}(q({\mathbf x})))$
for ${\mathbf x} \in (V_1 \otimes {\mathbb R})^n$ (the Gaussian vector) the
intgral $I(\tilde{g};\varphi_1)$ is the adelic function corresponding to the
Siegel modular form
 $$\vartheta^{(n)}~({\rm gen}(L_j),Z)$$
in the usual way.
 \vspace{0.3cm}\\
For the space $V_2$ we consider two different test functions at infinity: If
we choose a fixed majorant $\xi$ of $q$ and put 
 $$\varphi_{2,\infty,\xi}({\mathbf x}) = \exp(-2\pi~{\rm tr}(\xi({\mathbf x})))
 \quad \mbox{for } {\mathbf x} \in (V_2 \otimes {\mathbb R})^n,$$
the value of the theta kernel
 $$\theta(\tilde{g},{\mathbf 1}_{V_2},\varphi_{2,\infty,\xi} \otimes 
  \prod_{p \not= \infty} \varphi_{2,p})$$
at $h_2 = {\mathbf 1}_{V_2}$ corresponds to the theta function
 $$\vartheta^{(n)} (L_2,\xi,Z) = \sum_{{\mathbf x} \in L_2^n} \exp (2\pi i~
 {\rm tr}(q({\mathbf x})X)) ~\exp(-2 \pi~{\rm tr}(\xi({\mathbf x}) Y))$$
(with $Z = X+iY \in \mathfrak H_n$) considered by Siegel, and its integral
over $O_{(V_{2,q})}({\mathbb Q}) \setminus O_{(V_{2,q})}({\mathbb A})$
corresponds to the integral of this theta function over the space of
majorants $\xi$;
this is a nonholomorphic modular form in the space of Eisenstein
series by Siegel's theorem (or its extension to the Siegel-Weil-Theorem).
 \vspace{0.3cm}\\
Applying a certain differential operator as outlined in
\cite{kudla-compositio03} 
to $\varphi_{2,\infty,\xi}$, we obtain a different test function 
$\varphi'_{2,\infty,\xi}$, and the integral of the theta kernel
$\theta(\tilde{g},h,\varphi'_{2,\infty,\xi} \otimes \prod_{p \not= \infty} 
\varphi_{2,p})$ over $O_{(V_{2,q})}({\mathbb Q}) \setminus O_{(V_{2,q})}
({\mathbb A})$ corresponds to the holomorphic theta series of the indefinite
lattice $L_2$ considered by Siegel in \cite{Si} and by Maa\ss\ in \cite{Ma}
whenever the latter is defined.
 \vspace{0.3cm}\\
To simplify the discussion, we restrict now (following \cite{kudla-compositio03})
to $n = 1$. We denote by $\chi$ the quadratic character of ${\mathbb Q}_{\mathbb A}^
{\times}/{\mathbb Q}^{\times}$ defined by 
 $$\chi_v(x) = (x,(-1)^{\frac{m(m-1)}{2}} d)_v$$
for all places $v$, where $(~,~)_v$ is the Hilbert symbol. Then associated
to $\varphi_j$ there is a unique standard section $\Phi_j:\:
\tilde{G}({\mathbb A}) 
\times {\mathbb C} \longrightarrow {\mathbb C}$ with $\Phi_j(\cdot,s) \in
I(s,\chi)$, (where $I(s,\chi)$ is the principal series representation of 
$\tilde{G}({\mathbb A})$ with parameter $s$ and character $\chi$) such that
for $s_0 = \frac{m}{2}-1$ one has 
 $$\Phi_j(\tilde{g},s_0) = (\omega_j(\tilde{g})\varphi_j)(0) =: \lambda_j
  (\varphi_j).$$
With the Eisenstein series 
 $$E(\tilde{g},s;\varphi_j) := \sum_{\gamma \in \tilde{P}_{\mathbb Q} \setminus
  \tilde{G}_{\mathbb Q}} \Phi_j(\gamma\tilde{g},s)$$
associated to $\Phi_j$, the Siegel-Weil theorem asserts that 
$E(\tilde{g},s;\varphi_j)$ is holomorphic at $s = s_0$ and that one has the
identities 
 $$E(\tilde{g},s_0;\varphi_j) = \kappa \cdot I(\tilde{g};\varphi_j)$$
where $\kappa = 2$ if $m \leq 2$ and $\kappa = 1$ otherwise.
 \vspace{0.3cm}\\
The above maps $\lambda_j:\: S(V({\mathbb A})) \longrightarrow I(s_0,\chi)$
factor into a product $\lambda_j = \prod_{v} \lambda_{j,v}$ over all
places $v$ of ${\mathbb Q}$ and Kudla gives the following definition.

 \begin{kspdefinition} (Kudla)
 \begin{itemize}
 \item[(a)] Let $v$ be a (finite or infinite) place of ${\mathbb Q}$, let
$V_{1,v}$ and $V_{2,v}$ be quadratic spaces over ${\mathbb Q}_v$ of dimension
$m$ and discriminant $d$. Then functions $\varphi_{1v} \in S(V_{1,v})$
and $\varphi_{2,v} \in S(V_{2,v})$ are said to match if 
$\lambda_{1,v}(\varphi_{1,v}) = \lambda_{2,v} (\varphi_{2,v})$.
 \item[(b)] Let $V_1,V_2$ be quadratic spaces over ${\mathbb Q}$ of
the same dimension $m$ and discriminant $d$. Then two test functions
$\varphi_1 \in S(V_1({\mathbb A}))$ and $\varphi_2 \in S(V_2({\mathbb A}))$
match, if $\lambda_1(\varphi_1) = \lambda_2(\varphi_2)$. Equivalently, two
factorisable test functions $\varphi_1 = \bigotimes_{v} \varphi_{1,v}$,
$\varphi_2 = \bigotimes_{v} \varphi_{2,v}$ match if $\varphi_{1,v}$ and
$\varphi_{2,v}$ match for all places $v$.
 \end{itemize}
 \end{kspdefinition}
The matching principle observed by Kudla in \cite{kudla-compositio03} then 
states that for matching test functions $\varphi_1 \in S(V_1({\mathbb A}))$,
$\varphi_2 \in S(V_2({\mathbb A}))$ one has with $\Phi(\cdot,s_0) = \lambda_1(\varphi_1)
= \lambda_2(\varphi_2)$:
 $$I(\tilde{g};\varphi_1) = E(\tilde{g},s_0,\Phi) = I(\tilde{g};\varphi_2).$$
Although this identity is a trivial corollary of the Siegel-Weil
theorem, the matching 
principle gives highly nontrivial arithmetical identities since the integrals
$I(\tilde{g},\varphi_1)$ and $I(\tilde{g},\varphi_2)$ carry completely
different arithmetic information; in \cite{kudla-compositio03} the
principle is exploited 
to give identities between degrees of certain special cycles on
modular varieties 
and linear combinations of representation numbers of positive definite
quadratic forms. Kudla gives in \cite{kudla-compositio03} explicit
local  matching functions at the infinite place
and asserts the existence of local matching functions at the finite places
for $m > 4$ and for $m = 4$ if $\chi_p \not= 1$.
 \vspace{0.3cm}\\
We can now state the contribution of our computations from the
previous sections to this matching principle:
 \begin{kspproposition}
Let $L,V,q$ be as in the previous sections, let $n=1$ and let $\varphi_1 = \prod_{v}
\varphi_{1,v} \in S(V({\mathbb A}))$ be the test function for the positive
definite lattice $L$ as described above. Assume that $L$ is of square free odd 
level $N$ and that all $p|N$ divide the discriminant of $L$ to an odd power.
Let $\chi$ be the (primitive) quadratic character mod $N$ with 
$\vartheta(L,q) \in M_k(\Gamma_0(N),\chi)$ and let $p$ be a prime with $\chi(p)
= -1$. 

Let $\vartheta({\rm gen}(L))|_{T(p)} = \sum c_i \vartheta({\rm gen}(L_i))$
be the explicit linear combination of theta series of all the positive definite
genera of lattices of level $N$ and discriminant in $d \cdot ({\mathbb Q}^{\times})^2$
given by the results of Section 5, let $\psi_i$ be the test function attached
to the positive definite lattice $L_i$ as above. Let $(V_2,q_2)$ be the 
quadratic space $\tilde{V}$ of signature $(m-2,2)$ from Lemma 4.2 in Section 4,
let $L_2 = \tilde{L}$ in the notation of Lemma 4.2 and let
 $$\varphi_2' = \varphi'_{2,\infty,\xi} \otimes \prod_{p \not= 0}
 \varphi_{2,p}$$
be the test function attached to $L_2$ as described above.
 \vspace{0.3cm}\\
Then the test functions
 $$\psi := \sum_{i} c_i \psi_i \in S(V_1({\mathbb A}))$$
and 
 $$\varphi'_2 \in S(V_2({\mathbb A}))$$
match and we have
 $$I(\tilde{g},\psi) = I(\tilde{g},\varphi'_2).$$
 \end{kspproposition}
 \begin{proof}
This is clear from the discussion above and Theorem 4.3.
 \end{proof} 
 \begin{ksprem}
As already stated in \cite{kudla-compositio03} the matching principle
can easily be generalized 
to arbitrary $\widetilde{Sp}_n$. In the range of our results in Sections 4
and 5 we have then examples for the matching principle for general $n$ in 
the same way as described above.
 \end{ksprem}

\vskip1cm 
Hidenori Katsurada\\
Muroran Institute of Technology\\
27-1 Mizumoto\\
Muroran, 050-8585\\
Japan\\
{\tt hidenori@mmm.muroran-it.ac.jp}

\vskip0.5cm
Rainer Schulze-Pillot\\
Fachrichtung 6.1 Mathematik\\
Universit\"at des Saarlandes \\
Postfach 151150\\
66041 Saarbr\"ucken\\
Germany\\
{\tt schulzep@math.uni-sb.de}

\end{document}